\documentclass[12pt]{article}
\newcommand{\qed}{\ifmmode$\Box$\else{\unskip\nobreak\hfil
\penalty50\hskip1em\null\nobreak\hfil$\Box$
\parfillskip=0pt\finalhyphendemerits=0\endgraf}\fi}
\begin{document}
\renewcommand{\theequation}{\thesection.\arabic{equation}}
\newtheorem{Pa}{Paper}[section]
\newtheorem{Tm}[Pa]{{\bf Theorem}}
\newtheorem{La}[Pa]{{\bf Lemma}}
\newtheorem{Cy}[Pa]{{\bf Corollary}}
\newtheorem{Rk}[Pa]{{\bf Remark}}
\newtheorem{Pn}[Pa]{{\bf Proposition}}
\newtheorem{Ee}[Pa]{{\bf Example}}
\newtheorem{Dn}[Pa]{{\bf Definition}}
\newtheorem{I}[Pa]{{\bf}}
\newtheorem{Py}[Pa]{{\bf Property}}
\oddsidemargin 0in
\topmargin -0.5in

\begin{center}

{\bf \large Topologically free partial actions and faithful representations of partial crossed products }\\

\bigskip

 A.V. Lebedev\\

\bigskip

Belarus State University / University of Bialystok\\

\end{center}

\vspace{5mm}

\quad  \parbox{14,5cm}{\small \hspace{0.5cm}
In this paper we investigate the interrelation between the topological freedom of partial 
actions of discrete groups and 
faithful representations of partial crossed products

\medskip


\medskip

}
\vspace{5mm}
\tableofcontents
\section{Introduction}
\label{Intro}
\setcounter{equation}{0}

The notion of a partial crossed product of a $C^*$-algebra by an action of the  group $\bf Z$ 
by partial automorphisms was introduced by R. Exel \cite{Exel1}. It was then generalized further by 
K. McClanachan \cite{McClanachan} up to  partial crossed products by partial 
actions of discrete groups and by 
N. Sieben \cite{Sieben}   up to  partial crossed products by actions of inverse semigroups. 
A fruitful discussion of these and related objects one can  also find  in \cite{Paterson}. 
Partial crossed product is a natural generalization of the crossed product of a $C^*$-algebra by a group of automorphisms. To investigate this universal object it is important to have its faithful representations.
The description of the characteristic properties of such representations is the theme of this article.
Among the main properties in the presence of which one can obtain these 
representations is the existence  
of a contractive conditional expectation onto the 'coefficient' algebra 
(in this paper we call this property --- property $(^*)$)
and the topological freedom of the partial action. It is shown that the toplogical freedom implies property 
$(^*)$ and therefore gives us a powerful instrument to construct faithful representations of 
partial crossed products.\\

In this introductory section we gather the known necessary notions and facts on the partial actions and partial crossed products. In the next Section \ref{*top-free} we introduce the notion of the topologically free partial action and 
prove the principle result of the article (Theorem \ref{2.3}) linking toplogical freedom of the action and  property  $(^*)$.
Finally in Section \ref{3} on the base of this result we describe 
the existence of faithful representations of partial crossed products and reduced partial crossed products.\\

Let $A$ be  a $C^*$-algebra and $G$ be a discrete group.
We recall the definition of a partial $C^*$-dynamical systems and the corresponding 
partial crossed products for discrete groups (see, for example, \cite{McClanachan}).\\

{\bf Definition.} A {\em partial action of \/ $G$ on $A$} 
 (denoted by $\alpha $) is a collection $\{ D_g \}_{g\in G}$ of closed two-sided 
ideals of $A$ and a collection $\{ \alpha_g \}_{g\in G}$ of isomorphisms 
$\alpha_g : D_{g^{-1}} \to D_g $ such that 

$(1)\ \  \alpha_g \left( D_{g^{-1}} \cap D_h \right)\subset D_{gh} \ \ g,h \in G$

$(2) \ \ \alpha_{hg}(d) = \alpha_h (\alpha_g (d)) \ \ d\in D_{g^{-1}} \cap D_{g^{-1}h^{-1}}$

$(3)\ \  D_e =A, \ \ \alpha _e = {\rm Id}_A$

We shall say that $(A,G,\alpha)$ is a {\em partial dynamical system}.

Let 
$$
L = \left\{ a \in l^1 (G,A): a(g) \in D_g  \right\} .
$$
with the usual norm $\Vert a \Vert_1 = \sum\Vert a(g) \Vert$. 
Define a convolution  multiplication and involution on $L$ as follows.
$$
(a* b)(g)= \sum_{h\in G}\alpha_h \left[ \alpha_{h^{-1}}(a(h))b(h^{-1}g)   \right]
$$
$$
a^* (g)= \alpha_g (a(g^{-1})^*)
$$
With these operations $L$ becomes a Banach $^*$-algebra.

{\bf Definition.} \ \ The partial crossed product of $A$ and $G$
is the universal enveloping $C^*$-algebra of $L$. We 
denote the partial crossed product by $A\times_\alpha G$. 

{\bf Definition.} \ \  A {\em covariant representation} of $(A,G,\alpha)$ is a triple 
$(\pi ,u, H)$ where $\pi: A \to B(H)$ is a representation of $A$ on a Hilbert space $H$
 (here $B(H)$ is the algebra of all linear bouded operators on $H$),
\/ $u: G \to B(H)$ is a function $g \mapsto u_g$ with $u_g$ being 
a partial isometry on $H$ with the initial subspace $\left[  \pi (D_{g^{-1}}H) \right]$
and the final subspace $\left[  \pi (D_{g}H) \right]$ such that

$(1)\ \ u_g \pi (d)u_{g^{-1}} = \pi (\alpha_g(d))\ \ d\in D_{g^{-1}}$

$(2) \ \ \pi (d)\left[  u_gu_h -u_{gh}  \right]=0 \ \ d\in D_g\cap D_{gh}$

$3\ \ u^*_g = u_{g^{-1}}$.

{\bf Definition.} \ \ Let $(\pi ,u, H)$ be a covariant 
representation of $(A,G,\alpha)$. We define the  representation 
$\pi \times u : L \to B(H)$ by 
$$
(\pi\times u)(a) = \sum \pi (a (g))u_g .
$$
By the definition of $A\times_\alpha G$ \ \   $(\pi\times u)$ extends up to a $^*$-representation of 
$A\times_\alpha G$.
\begin{I}
\label{reduced}
\em
{\bf Reduced patrtial crossed product.}\ \ 
A special and important particular covariant representation of $(A,G,\alpha)$ is the so-called 
{\em reduced partial crossed product} which is defined in the following way (see \cite{McClanachan}, Section 3).

First we associate with  any representation $\pi : A \to B(H)$  a certain representation 
$
\tilde \pi 
$
 (the regular representation)  of $A$ on $l^2 (G,H)$. Let 
$$
\pi_g : D_g \to B(H)
$$ 
be defined by 
$$
\pi_g (d) =  \pi (\alpha _{g^{-1}} (d)).
$$
By \cite{Dixmier}, 2.10.4. there exists a unique extension $\pi_g^\prime$ of $\pi_g$ to $A$ which annihilates 
$[\pi_g (D_g)H]^\perp$. This extension is given by 
$$
\pi_g^\prime (d) = s \, \lim_\lambda \, \pi_g (v_\lambda d)
$$
where $\{ v_\lambda   \}_\lambda$ is an approximate identity for $D_g$. Now we define 
$$
\tilde \pi : A \to B (l^2 (G,H))
$$
by 
$$
\tilde \pi (d) \xi (g) = \pi_g^\prime (d) \xi (g) \ \ \  \xi \in l^2 (G,H), \ \ g\in G, \ \ d\in A.
$$
For the regular representation $\lambda$ of $G$  $\lambda : G \to B (l^2 (G,H))$ \ \ 
$(\lambda_g \xi ) (h) = \xi (g^{-1}h)$ we have (\cite{McClanachan} , Proposition 3.1.)  
$$
\lambda_g \tilde{\pi}(d) \lambda_{g^{-1}} = \tilde{\pi} (\alpha _g (d)) \ \ \ \ g\in G, \ \ d \in D_{g^{-1}}.
$$
If we let $ \tilde{\lambda_g} =  \lambda_gP_{g^{-1}}$ where $P_g$ is the orthogonal projection onto the Hilbert space
$[\tilde{\pi}(D_g)l^2 (G, H)]$  then $([\tilde{\pi}, \tilde{\lambda}, l^2 (G, H))$ is a covariant representation of 
$(A,G,\alpha )$.

Let $\Vert \cdot \Vert_r$ be the norm on $l^1 (G,A)$ defined by 
$$
\Vert a \Vert_r = \sup \{ \Vert  (\tilde{\pi}\times \lambda )(a)  \Vert (\pi , H)\in \/ {\rm Rep}\, (A) \}
$$
where ${\rm Rep}\, (A)$ is the set of all representations of $A$.

The  reduced partial crossed product $A\times_{\alpha r}G$ of $A$ by $G$ is the completion of $l^1 (G,A)$ 
with respect  to the norm $\Vert \cdot \Vert_r$.\\

In fact there is no need to use all representations of $A$ to define the reduced crossed product. As the next result tells it is enough to exploit any its faithful one.
\end{I}
\begin{Tm}
\label{1-} 
{\rm (\cite{McClanachan},  Proposition 3.4.)}\ \
Let $\pi : A \to B(H)$ be a representation. Then $\tilde{\pi}$ is faithful iff $\tilde{\pi}\times \lambda$ 
is faithful on $A\times_{\alpha r}G$.
\end{Tm}

\section{Property $(^*)$ and topologically free action}
\label{*top-free}
\setcounter{equation}{0}

{\bf Definition.} \ \ Let $(\pi ,u, H)$ be a covariant 
representation of $(A,G,\alpha)$. We shall say that $(\pi\times u)$ possesses 
{\em property } $(^*)$ if  for any finite sum
$$
\sum_{g\in F} \pi (a (g))u_g , \ \ F\subset G, \ \ \vert F\vert < \infty 
$$
we have 
$$
\left\Vert \sum_{g\in F} \pi (a (g))u_g  \right\Vert \ge \Vert a (e)\Vert 
$$
\begin{Rk}
\label{Rk-1}
\em It follows from  \cite{McClanachan},  Proposition 3.5 that 
$A\times_{\alpha r}G$ and  $A\times_\alpha G$ possess 
property $( ^*)$.
\end{Rk}
\begin{I}
\label{N}
\em
If $(\pi\times u)$ possesses property  $(^*)$ then the mapping 
$$
{\cal N}\left(\sum_{g\in F} \pi (a (g))u_g\right) = a(e) 
$$
is uniquely extended up to the mapping (positive, contractive, conditional expectation)
$$
{\cal N} : (\pi\times u) \to A .
$$
\end{I}
\begin{Rk}
\label{*-Exel}
\em
R. Exel \cite{Exel3}, Theorem 3.3.  proved that one can formulate property $(^*)$ in somewhat weaker 
but anyway equivalent way. In fact he proved a more general statement concerning graded $C^*$-algebras.
 Hereafter we formulate his result (its simplification) in terms of the objects considered in this paper.

{\em Let $(\pi\times u)$ be such that $\pi$ is a faithful representation of $A$ and 
$$
E : (\pi\times u) \to A
$$
be a bounded linear map such that 

a) $E (\pi (a)I) = a, \ \  a\in A, $

b) $E (\pi (a (g))u_g) = 0, \ \ g\neq 0 .   $

\noindent Then $(\pi\times u)$ possesses property $(^*)$ and $E={\cal N}$ where ${\cal N}$ is that mentioned in 
\ref{N}.
}
\end{Rk}
\begin{I}
\label{2.1} \em 
Now we proceed to  one of the main notions of the  article: topologically free action. To start with we note that  partial action defines in a natural way a partial dynamical system  (the action of a group by partial homeomorphisms) on the primitive ideal space 
${\rm Prim}\, A$ and the spectrum $\hat{A}$ of $A$. Here we give the description of this partial dynamical system.\\

For any ideal \/ $J\subset A$ we set  \ \ ${\rm supp}\, J = \{ x\in {\rm Prim}\, A : x  \not\supset  J \}$. 
It is known (see \cite{Dixmier}, 3.2.1.)  that  
the  mapping\ \  $x \to x\cap J$\ \ establishes a homeomorphism \ \ ${\rm supp}\, J \leftrightarrow {\rm Prim }\, J$ \ \ (with respect to the Jacobson topology) and \/ ${\rm supp}\, J$ is an open set in \/ ${\rm Prim }\, A$. Set also \/ $\hat{A}^J = \{ \pi \in \hat{A} : \pi (J) \neq 0     \}$ 
(here $\hat{A}$ is the spectrum of $A$). Then the mapping \/ $\pi \to \left. \pi  \right|_J$ establishes a homeomorphism 
$\hat{A}^J\leftrightarrow \hat{J}$ (with respect to the Jacobson topology) and $\hat{A}^J$ is an open set in $\hat{A}$  (see  \cite{Dixmier}, 3.2.1.).\\

Let us define the mapping \/ $\tau_g : \hat{A}^{D_{g^{-1}}} \to \hat{A}^{D_{g}}$ in the following way: 
for any \/ $\pi \in \hat{A}^{D_{g^{-1}}}$ we set 
$$
\tau_g (\pi ) (j) = \pi (\alpha_g^{-1}(j)), \ \ \ j\in D_g .
$$
The foregoing observations tell us that  $\tau_g$ is a homeomorphism.\\

Let us also define the mapping\  \/$t_g : {\rm supp}\, D_{g^{-1}} \to {\rm supp}\, D_g$ in the following way:
for any point \/ $x \in {\rm supp}\, D_{g^{-1}}$ such that \/ $x = {\rm ker}\, \pi$ where $\pi \in  \hat{A}^{D_{g^{-1}}}$ 
we set 
$$
t_g (x) = {\rm ker }\, \tau_g (\pi ).
$$
Clearly \/$t_g$\/ is a homeomorphism.\\

For $\tau_g$ and $t_g$ defined in the above deiscribed way we have that $\{ \tau_g  \}_{g\in G}$ defines an action of $G$  by  {\em partial homeomorphims} of $\hat{A}$ and $\{  t_g \}_{g\in G}$  defines an action of $G$  by  {\em partial homeomorphisms} of 
${\rm Prim}\, A$.
\end{I}
\begin{I}
\label{top-free}
\em 
We say that the action $\{ \alpha _g \}_{g \in G}$ is {\em toplogically free} iff for any finite set 
$\{ g_1 , ... g_k \} \subset G$ and any nonempty open set
 \/ $U \subset {\rm supp}\, D_{g_1^{-1}}\cap ... \cap   {\rm supp}\, D_{g_k^{-1}}$ \/  
there exists a point $x \in U $ such that all the points $t_{g_i}(x), \ \  i=\overline{1,k}$ are distinct.

This condition can be also formulated in the following way: for any finite set 
$\{ g_1 , ... g_k \} \subset G$ and any nonempty open set $U$ there exists a point $x \in U $ such that all the points $t_{g_i}(x), \ \  i=\overline{1,k}$\ \  that  are defined  ($ \Leftrightarrow x \in {\rm supp}\, D_{g_i^{-1}}$) are distinct.

If we denote by $X_g$ the set 
$$
X_g = \{x \in {\rm supp}\, D_{g^{-1}}: t_g (x)=x   \}
$$
then the foregoing condition can be also written in the next  way:
 for any finite set $\{ g_1, ... g_n \}, \ \ g_i \neq e $ the interiour of the set  
$
\left[  \cup_{i=1}^n X_{g_i}
  \right]
$
is empty. 
\end{I}
The main statement of this section is Theorem \ref{2.3} and the most  important technical result  is Lemma \ref{2.6}.   Among the technical instruments of the proof 
of this lemma  is the next Lemma \ref{a} which is useful in its own right.
\begin{La}
\label{a}
{\rm (\cite{AntLeb}, Lemma 12.15).} Let $B$ be a $C^*$-subalgebra of the algebra $L(H)$ of linear bounded operators in a Hilbert space $H$. If \/ $P_1, \  P_2\  \in\  B^\prime$ are two orthogonal projections such that the restrictions 
$$
\left. B \right|_{H_{P_1}} \ \ and \ \  \left. B \right|_{H_{P_2}}
$$
(where $H_{P_1} = P_1 (H), \ \ H_{P_2} = P_2 (H)$) are both irreducible and these restrictions are distinct representations then 
$$
H_{P_1} \perp  H_{P_2} .
$$
\end{La}
\begin{La}
\label{2.6}
Let the action  $\{ \alpha _g \}_{g \in G}$ be topologically free and  $(\pi\times u)$  is such that $\pi$ is a faithful representation of $A$.  Let   $F $ be a finite subset of $G$, and  
  $a \in L$  be  any function such that $a(g)=0, \ \ g\notin F$, and $c\in (\pi\times u)$ be the  operator   of the form 
\begin{equation}
\label{e--}
c= \sum_{g\in F} \pi (a (g))u_g  .
\end{equation}
Then for every $\varepsilon >0$ there exists an irreducible representation $\pi^\prime$ of $\pi (A)$ such that for any  irreducible 
 representation $\nu$ of $(\pi\times u)$ which is an extension of $\pi^\prime$ we have

(i)\ \ $\Vert \pi^\prime [\pi (a (e))]\Vert \ge \Vert a (e) \Vert  - \varepsilon$,

(ii)\ \ $  P_{\pi^\prime } \,   \pi^\prime [\pi (a (e))]\, P_{\pi^\prime } = P_{\pi^\prime } \, \nu (c)\,  P_{\pi^\prime }  $

\noindent where $P_{\pi^\prime }$ is the orthogonal projection onto $H_{\pi^\prime }$ in $H_\nu$.
\end{La}
{\bf Proof.}\ \  As $\pi (A) \cong  A$ we shall identify throughout the proof $\pi (A)$ and $A$ (in order to shorten the notation). 

For any $d\in A$ and $x\in {\rm Prim}\, A$ we denote by $\breve{d} (x) $ the number 
\begin{equation}
\label{e2.1}
\breve{d} (x)  = \inf_{j\in x}\Vert  d+j \Vert
\end{equation}
For every $d\in A$ the function $\breve{d} (x) $ is lower semicontinuous on ${\rm Prim}\, A$ and attains its upper bound equal to $\Vert d \Vert$ (see \cite{Dixmier}, 3.3.2. and 3.3.6.). 

Let $x_0 \in {\rm Prim}\, A$ be a point at which  $\breve{a}(e) (x_0) =  \Vert a(e) \Vert$ and $\pi_0$ be an irreducible representation of $A$ such that $x_0 = {\rm ker}\, \pi_0$ (thus $\Vert \pi_0 (a(e))\Vert = \Vert a(e) \Vert$).
 Since the function $\breve{a}(e) (x)$ is lower semicontinuous it follows that for any $\varepsilon >0$ 
there exists an open set 
$U \subset {\rm Prim}\, A$ such that 
\begin{equation}
\label{e2.1a}
\breve{a}(e) (x) >  \Vert a(e) \Vert - \varepsilon \ \ {\rm for\ \  every}\ \ x\in U. 
\end{equation}
As 
 the action  $\{ \alpha _g \}_{g \in G}$ is topologically free there exists a point 
$x^\prime \in U $ such that all the points $t_{g_i}(x^\prime ), \ \  i=\overline{1,k}$ are   distinct (if they are defined
$\Leftrightarrow x^\prime \notin {\rm supp}\, D_{g_i^{-1}}$).

Let $\pi^\prime$ be an 
irreducible representation of $A$ such that ${\rm ker}\, \pi^\prime = x^\prime$ and let $\nu$ be any extension of 
$\pi^\prime$ up to an irreducible representation of   $(\pi\times u)(L)$. We shall denote by the same letter $\nu$ an extension of the mentioned representation up to an irreducible representation of the $C^*$-algebra  $C$ generated by $(\pi\times u)(L)$ and $\{u_g \}_{g\in G}$ (see \cite{Dixmier}, 2.10.2.). 
 For this representation $\nu$ we have
$$
H_{\pi^\prime}\subset H_\nu
$$
where $H_{\pi^\prime}$ is the representation space for $\pi^\prime$ and $H_\nu$ is that for $\nu$.

By the choice of $\pi^\prime$ and (\ref{e2.1a}) we conclude that there exists a vector  $\xi \in H_{\pi^\prime}$ 
such that $\Vert \xi \Vert =1$ and 
\begin{equation}
\label{e00}
\Vert {\pi^\prime} (a(e))\xi \Vert  > \Vert a(e)\Vert -\varepsilon.
\end{equation}
Thus (i) is proved.

To prove (ii)  let us  observe first that for any vectors $\xi , \eta \in  H_{\pi^\prime}$ we have 
\begin{equation}
\label{e2.2}
\left\langle {\pi^\prime}(d_1)\eta ,\,  \nu (d_2 u_g)\xi \right\rangle =0, \ \ \ \ d_1 \in A,\ \ d_2 \in D_g, \ \ g\in F,\ \ g\neq e .
\end{equation}
Which in turn will imply 
\begin{equation}
\label{e2.2a}
P_{\pi^\prime}\,  \nu (d_2 u_g)\, P_{\pi^\prime} =0, \ \  g\in F, \ \ d_2 \in D_g, \ \  g\neq e
\end{equation}

To prove (\ref{e2.2})  we consider the following possible positions of $x^\prime$.\\

$x^\prime \notin {\rm supp}\, D_g$. \ \  In this case we have ${\pi^\prime}(d_2^*)=0$ and 
$$
\left\langle {\pi^\prime}(d_1)\eta ,\,  \nu (d_2) \nu (u_g)\xi \right\rangle =
\left\langle \nu (d_2^*){\pi^\prime}(d_1)\eta ,\,   \nu (u_g)\xi \right\rangle =
$$
$$
\left\langle \pi^\prime (d_2^*){\pi^\prime}(d_1)\eta ,\,   \nu (u_g)\xi \right\rangle = 0.
$$
\ \\

$x^\prime \notin {\rm supp}\, D_{g^{-1}}$. \ \  Observing that $\nu (u_g^*u_g)$ 
is the projection onto the essential space of $\nu (D_{g^{-1}})$ we conclude that 
$\nu (u_g^*u_g) \xi =0$ and therefore we have 
$$
\left\langle {\pi^\prime}(d_1)\eta ,\,  \nu (d_2) \nu (u_g)\xi \right\rangle =
\left\langle {\pi^\prime}(d_1)\eta ,\,  \nu (d_2) \nu (u_gu_g^* u_g)\xi \right\rangle =
$$
$$
\left\langle {\pi^\prime}(d_1)\eta ,\,  \nu (d_2) \nu (u_g) \nu  (u_g^*u_g)\xi \right\rangle = 0 .
$$
\ \\

Finally let $x^\prime \in \left[{\rm supp}\, D_g \cap {\rm supp}\, D_{g^{-1}}\right]$.\\
 In this case $\pi^\prime$ 
is an irreducible representation as for $D_g$ so also for $D_{g^{-1}}$ and 
$t_g (x^\prime )\in {\rm supp}\, D_g$ (according to the definition of $t_g$ \ref{2.1}). Moreover  we have 
\begin{equation}
\label{e2.3}
\nu (u_g^* u_g)\eta = \eta , \ \ \  \nu (u_g u_g^* )\eta = \eta  \ \ {\rm for \ \ any}\ \ \eta \in H_{\pi^\prime}.
\end{equation}
Since $\nu (u_g)$ is a partial isometry the observation (\ref{e2.3}) implies that $H_{\pi^\prime}$ belongs as to the initial and final subspaces of $\nu (u_g)$ so also to the initial and final subspaces of $\nu (u_g^*)$ and the mappings 
\begin{equation}
\label{e2.4}
\nu (u_g) : H_{\pi^\prime} \to \nu (u_g) \left[ H_{\pi^\prime}\right] \ \ {\rm and}\ \ 
\nu (u_g^*) : H_{\pi^\prime} \to \nu (u_g^*) \left[ H_{\pi^\prime}\right] 
\end{equation}
are isomorphisms.

Let $P_1$ be the orthogonal projection of $H_\nu$ onto $H_{\pi^\prime}$. By the definition of $\nu$ we have that 
$P_1\in \nu (A)^\prime$  and (\ref{e2.3}) means that 
\begin{equation}
\label{e2.5}
P_1 = P_1\, \nu (u_g^* u_g) = P_1\,  \nu (u_g u_g^* ) .
\end{equation}

Set $P_2 = \nu (u_g)P_1\nu (u_g^*)$. The foregoing observations imply that 
$$
\nu (u_g): P_1(H_{\pi^\prime}) \to P_2(H_{\pi^\prime})
$$
 is an isomorphism. Observe also that 
\begin{equation}
\label{e2.6}
P_2 \in \left[ \nu (D_g)  \right]^\prime .
\end{equation}
Indeed. For any $d\in D_g$ we have
$$
\nu (d) = \nu  (u_g u_g^* ) \nu (d) = \nu (d)\nu  (u_g u_g^* )
$$
and 
$$
\nu (\alpha_{g^{-1}}(d)) = \nu (u_g^* u_g) \nu (\alpha_{g^{-1}}(d))= \nu (\alpha_{g^{-1}}(d)) \nu (u_g^* u_g),
$$
and 
$$
 \nu (u_g^*) \nu (d) \nu (u_g)= \nu (\alpha_{g^{-1}}(d)),
$$
and 
$$
\nu (\alpha_g \circ \alpha_{g^{-1}}(d)) =\nu (a).
$$
Using this we obtain for any $d\in D_g$
\begin{eqnarray*}
P_2 \nu (d)  &=\nu (u_g)P_1 \nu (u_g^*)\nu (d) =  \nu (u_g)P_1 \nu (u_g^*) \nu (u_g u_g^*)\nu (d) =\\
  &\nu (u_g)P_1 \left[\nu (u_g^*)\nu (d) \nu (u_g)\right] \nu (u_g^*) =  \nu (u_g)P_1 \nu (\alpha_{g^{-1}}(d))\nu (u_g^*)=\\
 & \nu (u_g)\nu (\alpha_{g^{-1}}(d))\nu (u_g^*u_g)P_1 \nu (u_g^*)= 
\nu (\alpha_g \circ\alpha_{g^{-1}}(d))\nu (u_g)P_1\nu (u_g^*)=\\ \nu (d)P_2 & \ \ \\
\end{eqnarray*}
Thus (\ref{e2.6}) is true.

In addition the irreducibility of 
$\left. {\nu (D_g)} \right|_{H_{P_1}}$ implies the irreducibility of 
$\left. {\nu (D_g)} \right|_{H_{P_2}}$ (here $H_{P_1} = P_1 (H_\nu )=H_{\pi^\prime}$ and 
$H_{P_2} = P_2 (H_\nu )$).

Now observe that for $d\in D_g$ we have 
$$
P_1\,  \nu (d) = 0  \Leftrightarrow \breve{d}(x^\prime )=0
$$
and 
\begin{eqnarray*}
P_2\,  \nu (d) =0 \Leftrightarrow  \nu (u_g)\, P_1\,  \nu (u_g^*)\nu (u_g)\nu(u_g^*) \nu (d) = 0
\Leftrightarrow \\
\nu (u_g)\nu (u_g^*u_g)\, P_1 \,  \nu \left( \alpha_{g^{-1}} (d)\right)\nu(u_g^*) = 0
\Leftrightarrow \\
\nu (u_g)\nu (u_g^*u_g)\, P_1 \,  \nu \left( \alpha_{g^{-1}} (d)\right)\, P_1\,  \nu (u_g^*u_g)\nu (u_g^*) = 0
\Leftrightarrow \\
P_1 \,  \nu \left( \alpha_{g^{-1}} (d)\right) =0 \Leftrightarrow 
\mathord{\buildrel{\lower3pt\hbox{$\scriptscriptstyle\smile$}} 
\over {\alpha_{g^{-1}}}(d)} (x^\prime)=0\Leftrightarrow \breve{d} (t_g (x^\prime)) =0.
\end{eqnarray*}
So (since the points $x^\prime$ and $t_g (x^\prime)$ are distinct) we conclude that the representations 
$\left. {\nu (D_g)} \right|_{H_{P_1}}$ and $\left. {\nu (D_g)} \right|_{H_{P_2}}$ are distinct. Applying Lemma \ref{a}
we find that 
\begin{equation}
\label{p1p2}
P_1\cdot P_2 =0.
\end{equation}
 By applying (\ref{p1p2}), (\ref{e2.3}), (\ref{e2.5}) and  (\ref{e2.6}) we have  for any 
$\eta, \xi \in H_{P_1} =H_{\pi^\prime}, \ \ d_1 \in A,\ \ d_2 \in D_g$ 
\begin{eqnarray*}
\left\langle {\pi^\prime}(d_1)\eta ,\,  \nu (d_2 u_g)\xi \right\rangle = 
\left\langle P_1\, \nu(d_1)\eta ,\,  \nu (d_2 u_g)\, P_1\,  \nu (u_g^*u_g)\xi \right\rangle = \\
\left\langle P_1\,\nu(d_1)\eta ,\,   \nu (d_2)\, P_2 \, \,  \nu (u_g) \xi \right\rangle =
\left\langle P_1\,\nu(d_1)\eta ,\,  P_2\, \, \nu (d_2) \nu (u_g) \xi \right\rangle =\\
\left\langle P_2\cdot P_1\,\nu(d_1)\eta ,\,   \nu (d_2) \nu (u_g) \xi \right\rangle = 0 
\end{eqnarray*}
which finishes the proof of  (\ref{e2.2}) (and therefore the proof of (\ref{e2.2a}) as well).

Now returning to  the operator (\ref{e--})  (recall that we are identifying $A$ and $\pi (A)$) 
and using (\ref{e2.2a}) we have that 
$$
P_{\pi^\prime}\,   \nu \left[\sum_{g\in F} a (g) u_g \right]\,  P_{\pi^\prime} = 
P_{\pi^\prime} \,  \nu (a (e))\,  P_{\pi^\prime} = P_{\pi^\prime} \,   {\pi^\prime}(a (e)) \, P_{\pi^\prime}
$$
so (ii) is true  and the proof of the lemma is complete.\\

% (\ref{e1})  and using (\ref{e00}) and (\ref{e2.2}) we have that 
%\begin{eqnarray*}
%\left\Vert \sum_{g\in F} \pi (a (g))u_g  \right\Vert ^2 \ge 
%\left\Vert \nu \left(\sum_{g\in F} \pi (a (g))u_g \right) \xi \right\Vert ^2 = \\
%\Vert  \nu (a(e))\xi \Vert^2 +  \left\Vert \nu \left(\sum_{g\neq e} \pi (a (g))u_g \right) \xi \right\Vert ^2 \ge \\
%\Vert  \nu (a(e))\xi \Vert^2  = \Vert  \pi^\prime (a(e))\xi \Vert^2 \ge \left(  \Vert a(e) \Vert - 2\varepsilon  \right)^2 .
%\end{eqnarray*}
%In view of the arbitrarness of $\varepsilon$ this finishes the proof of the theorem.

As an immediate corollary of this lemma we obtain the next 
\begin{Tm}
\label{2.3}
Let the action  $\{ \alpha _g \}_{g \in G}$ be topologically free. If $(\pi\times u)$  is such that $\pi$ is a faithful representation of $A$ then $(\pi\times u)$ possesses property $(^*)$.
\end{Tm}
{\bf Proof.}\ \ Let $c$ be the operator  (\ref{e--}). Take $\pi^\prime $ mentioned in 
 the statement of Lemma \ref{2.6}. 
Then we have by (ii) and (i) 
$$
\Vert c \Vert \ge \Vert \nu (c)\Vert \ge \Vert P_{\pi^\prime}\, \nu (c) \, P_{\pi^\prime} \Vert \ge 
\Vert a(e) \Vert - \varepsilon
$$
In view of the arbitrarness of $\varepsilon$ this implies property $(^*)$.\\

One more simple corollary of Lemma \ref{2.6} and Theorem \ref{2.3} is the following 
\begin{La} 
\label{2.6a}
Let the action  $\{ \alpha _g \}_{g \in G}$ be topologically free and  $(\pi\times u)$  is such that $\pi$ is a faithful representation of $A$.  
Then for any $c\in (\pi\times u)$ and every $\varepsilon >0$ there exists an irreducible representation $\pi^\prime$ of $\pi (A)$ such that for any  irreducible 
 representation $\nu$ of $(\pi\times u)$ which is an extension of $\pi^\prime$ we have

(i)\ \ $\Vert \pi^\prime [{\cal N}(c)]\Vert \ge \Vert {\cal N}(c) \Vert  - \varepsilon$,

(ii)\ \ $ \Vert  P_{\pi^\prime } \,   \pi^\prime \left[{\cal N}(c)\right]\, P_{\pi^\prime } - P_{\pi^\prime } \, \nu (c)\,  P_{\pi^\prime } \Vert \le \varepsilon $
\end{La}
{\bf Proof.}\ \  Follows from the standard approximation argument in view of the density of finite sums of the form 
(\ref{e--}) in $(\pi\times u)$ and the fact that $(\pi\times u)$ possesses property $(^*)$.

\section{Property $(^*)$, topologically free action, partial crossed products and partial reduced crossed products}
\label{3}
\setcounter{equation}{0}

It is reasonable to consider $A\times_\alpha G$ as the {\em maximal} $C^*$-algebra pssessing property $(^*)$
(it follows from the construction of  $A\times_\alpha G$ and Remark \ref{Rk-1}). On the other hand 
it has been shown by 
R. Exel that  $A\times_{\alpha r} G$ is the {\em minimal}\ \ $C^*$-algebra pssessing this property.
The exact meaning of  'minimality' is given in the next statement which is a reformulation 
(in fact simplification) of \cite{Exel3}, Theorem 3.3. (we recall at this point that according to 
Theorem \ref{1-}
 for any 
faithful representation $\pi$ of $A$\ \   $\tilde{\pi}\times \lambda$  is a faithful representation of 
$A\times_{\alpha r} G$).
\begin{Tm}
\label{min}
Let   $(\pi\times u)$  be such that $\pi$ is a faithful representation of $A$.  If  
$(\pi\times u)$ possesses property $(^*)$  then the mapping 
$$
(\pi\times u) \ni \sum \pi (a (g))u_g \mapsto \sum \tilde{\pi}(a(g))\tilde{\lambda}_g \in 
\tilde{\pi}\times \lambda
$$
can be extended up to a $C^*$-algebra epimorphism (here $\tilde{\pi}\times \lambda$ is that mentioned in Theorem \ref{1-}). 
\end{Tm}
\begin{Rk}
\label{cro-red}
\em
It is also known that if $G$ is an amenable group then the canonical surjection $\Lambda : A\times_\alpha G
\to A\times_{\alpha r} G$ is an isomorphism (see, for example, \cite{McClanachan}, Proposition 4.2).
\end{Rk}
This observation along with Theorem \ref{min} leads to the next result 
\begin{Tm}
\label{isom-*}
 Let  $G$ be an amenable group and
$(\pi^i ,u^i, H^i), \ \ i=1,2$ be two covariant representations of  $(A,G,\alpha)$  such that both 
$\pi^i \times u^i, \ \ i=1,2$ possess property $(^*)$ then the mapping 
$$
\sum \pi^1 (a (g))u^1_g \mapsto \sum \pi^2 (a (g))u^2_g 
$$
give rise to the isomorphism of the algebras $\pi^1 \times u^1$ and $\pi^2 \times u^2$.\\
\end{Tm}
\begin{Rk}
\label{istorija-*}
\em 
  The importance of  property $(*)$ for the first time (probably) was clarified
by O'Donovan \cite{O'Donov} in connection with the description of  $C^*-$al\-ge\-bras generated by
weighted shifts. The most general result  (of Theorem \ref{isom-*} type) establishing the crucial role of this property in the
theory of crossed products of $C^*-$algebras by discrete groups of {\em automorphisms} was
obtained in \cite{Leb1} for an {\em arbitrary} $C^*-$algebra and {\em amenable} discrete group (see also \cite{AntLeb}, Chapters 2,3\ \ for complete proofs and various
applications). The relation
of the corresponding property to the faithful representations of crossed products by {\em
endomorphisms} generated by {\em isometries} was investigated in \cite{BKR,ALNR}. 

It is worth mentioning that in \cite{AntLeb}, Theorem 12.8 (an analogue to  Theorem \ref{isom-*}) was proved 
in a direct way {\em not} exploiting the reduced crossed product  so in particular  the isomorphism  
 of $\Lambda : A\times_\alpha G
\to A\times_{\alpha r} G$ for amenable groups can also be derived  from this result  (the proof of \cite{AntLeb}, Theorem 12.8 can be easily extended up to a partial crossed product situation).\\
\end{Rk}

Theorem \ref{2.3} gives us a possibility to verify property $(^*)$\ \ in an   automatic way  by means of the property
of the underlying partial dynamical system. This theorem along with the foregoing results leads to the following 
Theorems  \ref{3.1}, \ref {3.2}.
\begin{Tm}
\label{3.1}
Let the action  $\{ \alpha _g \}_{g \in G}$ be topologically free and  $(\pi\times u)$  is such that $\pi$ is a faithful representation of $A$.   Then the mapping 
$$
(\pi\times u) \ni \sum \pi (a (g))u_g \mapsto \sum \tilde{\pi}(a(g))\tilde{\lambda}_g \in 
\tilde{\pi}\times \lambda
$$
can be extended up to a $C^*$-algebra epimorphism (here $\tilde{\pi}\times \lambda$ is that mentioned in Theorem \ref{1-}). 
\end{Tm}
{\bf Proof.}\ \ Apply  Theorem \ref{2.3} and Theorem \ref{min}.
\begin{Tm}
\label{3.2}
Let $G$ be an amenable group and the action  $\{ \alpha _g \}_{g \in G}$ be topologically free. 
If $(\pi^i ,u^i, H^i), \ \ i=1,2$ be two covariant representations of  $(A,G,\alpha)$  such that both 
$\pi^i , \ \ i=1,2$ are faithful representations of $A$ then the mapping 
$$
\sum \pi^1 (a (g))u^1_g \mapsto \sum \pi^2 (a (g))u^2_g 
$$
give rise to the isomorphism of the algebras $\pi^1 \times u^1$ and $\pi^2 \times u^2$.
\end{Tm}
{\bf Proof.}\ \ Apply  Theorem \ref{2.3} and Theorem \ref{isom-*}.\\

The next Theorem \ref{2.6b} and Corollary \ref{2.7}\ \  are in a way opposite to Theorem \ref{min}.
They form a generalization of 
   \cite{Exel2},  Theorem 2.6.
(where $A=C_0 (X)$).
\begin{Tm}
\label{2.6b}
Let the action  $\{ \alpha _g \}_{g \in G}$ be topologically free. If $I$ is an ideal in $A\times_{\alpha r}G$ then 
$I = \{ 0 \}$ iff $I\cap A = \{ 0 \}$.
\end{Tm}
{\bf Proof.}\ \  Let  $I\cap A = \{ 0 \}$. Denote by $\pi : A\times_{\alpha r}G \to (A\times_{\alpha r}G)/I$ the quotient map
and let $c\in I$ be an element such that $c\ge 0$ and $\pi (c)=0$. To prove that $I = \{ 0 \}$ we have to verify that
\begin{equation}
\label{e2.11}
c=0 .
\end{equation}
Since the mapping 
$$
{\cal N} : A\times_{\alpha r}G \to A
$$
defined in  \ref{N} is faithful (see, for example,  \cite{Exel3}, Proposition 2.12) (\ref{e2.11}) will be proved if we prove that 
\begin{equation}
\label{e2.12}
{\cal N}(c) =0 .
\end{equation}
So let us verify the latter property.

Since $I\cap A = \{ 0 \}$ it follows that $\pi (A) \cong A$.
Given $\varepsilon > 0$ take $\pi^\prime$ form the statement of Lemma \ref{2.6a}  
(we can reffer to this representation either as to the  repersentation of $\pi (A) $ 
so also as to the representation of $A$) and 
  extend it up to an irreducible representation 
$\nu$ of  $ \pi  \left(A\times_{\alpha r}G\right)$ 
(here we consider $ \pi \left(A\times_{\alpha r}G \right)$ as $(\pi\times u)$ in the statement of Lemma \ref{2.6a}).  
Evidentely $\nu\circ\pi$ is an irreducible representation of $A\times_{\alpha r}G$.

Now the condition $\pi (c)=0$ and property (ii) of the statement of Lemma \ref{2.6a} imply
$$
\varepsilon \ge  \Vert  P_{\pi^\prime } \,   \pi^\prime \left[{\cal N}(\pi (c))\right]\, P_{\pi^\prime } - P_{\pi^\prime } \, \nu (\pi (c))\,  P_{\pi^\prime } \Vert = \Vert \pi^\prime [{\cal N}(\pi (c))]\Vert = \Vert \pi^\prime [{\cal N}(c)]\Vert .
$$
This and  (i) implies 
$$
\Vert {\cal N}(c)  \Vert \le 2\varepsilon .
$$
Which proves (\ref{e2.12}) by the arbitrariness of $\varepsilon$.
\begin{Cy}
\label{2.7}
Let the action  $\{ \alpha _g \}_{g \in G}$ be topologically free. A representation $\pi$ of the reduced partial crossed product $A\times_{\alpha r}G$ is faithful iff it is faithful on $A$.
\end{Cy}
{\bf Proof.}\ \ Take in the statement of Theorem \ref{2.6b} $I = {\rm ker}\, \pi$. \\

\begin{Rk}
\label{top-free-istor}
\em
The interrelation between the topological freedom of the action and 
property $(^*)$ and application of these properties to various crossed product type results 
have been intensively exploited by many authors.
The treatment of the topological freedom as an instrument of 
investigation of ideals in the crossed products was started (probably) 
by D.P. O'Donovan in \cite{O'Donov}, Theorem 1.2.1. \ \ 
  Theorem \ref{3.2}
in the case of  a commutative algebra $A$ and the  action of the group $\bf Z$ by {\em automorphisms}
was proved in \cite{LebUMN,Leb-dis}. 
The development of this field and its numerious  (not purely $C^*$-algebraic) applications   such as, 
for example, the construction of  symbolic calculus and  the
solvability theory of functional differential equations one can find in \cite{Ant,AntLeb,AnLebBel1, AnLebBel2}.
For the general {\em automorphism} situation  Theorem \ref{3.2}   was
obtained in \cite{Leb1}  (see also in this connection \cite{AntLeb}, Chapters 2,3). 
Among the already mentioned 'purely' $C^*$-algebraic sources we have to emphasize an   outstanding 
  contribution to the theme made  in \cite{Exel3}. A deep and versatile study of the topological freedom
(in the situation $A=C_0 (X)$) and its  application to a series 
of structural problems in partial crossed product theory is implemented in \cite{Exel2}.

In the Lebesgue space situation the topological freedom corresponds to the so-called {\em metrical freedom}.
The interrelation between  this property, property $(^*)$ and the corresponding crossed product results 
(in the {\em automorphisms} situation)
were investigated and  applied to the solution 
of the  problem of classification of measure preserving automorphisms by W.B Arveson and K.B. Josephson in 
\cite{Arv,Arv-Jos}. 

In the {\em endomorphisms} situation namely in the case when a $C^*$-algebra endomorphism is generated by a
{\em single} isometry the interrelations between the topological freedom 
of the action and property $(^*)$ have been investigated  in \cite{Leb-UMN2,Leb1,Leb-end} where in particular the 
analogues to 
Theorems \ref{2.3}, \ref{isom-*}, \ref{3.2}    for the situation considered were obtained. In fact this research has been 
inspired by the pioneering work by  V.A. Arzumanian and  A.M. Vershik \cite{Arz-Ver1,Arz-dis,Arz-Ver2,Arz2}
where the corresponding Lebesgue space objects have been introduced and studied. 
Recently this theme has got a new development in the work by R. Exel \cite{Exel4,Exel5}, 
and R. Exel  and A.M. Vershik \cite{Exel-Ver}.
\end{Rk}

%??????????

%\begin{Pn}
%\label{*g}

%????????

%Let $(\pi\times u)$ possess  
% property  $( ^*)$ then for any finite sum 
% $$
% \sum_{g\in F} \pi (a (g))u_g , \ \ F\subset G, \ \ \vert F\vert < \infty 
 %$$
 %and any $g_0 \in G$ we have 
 %$$
 %\left\Vert \sum_{g\in F} \pi (a (g))u_g  \right\Vert \ge \Vert a (g_0)\Vert 
%$$
%\end{Pn}
%{\bf Proof.} \ \ 

%{\bf Definition.} \ \ The functional with finite support ???????

%\begin{Tm}
%\label{lin-funct}
%Let  $G$ be an amenable group then the set of 
%positive functionals with finite support is dense in the set of all 
%positive functional on $A\times_\alpha G$ 
%\end{Tm}
%{\bf Proof.}
%\begin{Tm}
%\label{isom}
%Let $G$ be an amenable group $\pi : A \to L(H)$ be a faithful representation of $A$ 
%and $(\pi ,u, H)$ be a covariant 
%representation of $(A,G,\alpha)$. Then $(\pi\times u)$ is a faithful representation of 
%$A\times_\alpha G$  iff \/ $(\pi\times u)$\/ possesses property $( ^*)$.
%\end{Tm}
%{\bf Proof.}\ \ 

%????????

\end{document}